\documentclass[letter,12pt]{article}

\newcounter{rmk}

\usepackage{comment}
\usepackage[linktocpage=true]{hyperref}
\usepackage[margin=1.2 in, top=1 in, bottom= 1.2 in]{geometry}

\usepackage[english]{babel}
\usepackage{amsfonts,amsmath}
\usepackage{amsthm}

\usepackage{mathrsfs}
\usepackage{bigints}
\usepackage[mathscr]{euscript}
\usepackage{bbm}
\usepackage{latexsym}
\usepackage{math dots}
\usepackage{amssymb}
\usepackage{mathtools}

\usepackage{enumitem}
\setlist{nolistsep}
\usepackage{amsthm}
\usepackage{relsize}
\usepackage{nicefrac}
\usepackage[capitalize]{cleveref}

\newtheoremstyle{plain}{3mm}{3mm}{\slshape}{}{\bfseries}{.}{.5em}{}
\newtheoremstyle{definition}{2mm}{2mm}{}{}{\bfseries}{.}{.5em}{}
\theoremstyle{plain}
	
\newtheorem{theorem}{Theorem}

\newtheorem{lemma}[theorem]{Lemma}

\theoremstyle{definition}

\newtheorem{remark}[rmk]{Remark}

\theoremstyle{plain}
\newtheorem*{namedthm}{\namedthmname}
\newcounter{namedthm}
	

\newcommand{\otherlim}[2]{\underset{#2}{{#1}\text{-}\lim}\ }
\newcommand{\supp}{{\normalfont\text{supp}}}
\newcommand{\Oh}{{\rm O}}
\newcommand{\oh}{{\rm o}}
\newcommand{\N}{\mathbb{N}}
\newcommand{\Z}{\mathbb{Z}}
\newcommand{\R}{\mathbb{R}}
\newcommand{\C}{\mathscr{C}}
\newcommand{\Q}{\mathbb{Q}}
\newcommand{\T}{\mathbb{T}}
\newcommand{\Hilb}{\mathscr{H}}
\newcommand{\Aut}{\mathsf{Aut}}
\newcommand{\Cont}{\mathsf{C}}
\newcommand{\Lip}{\mathsf{Lip}}
\newcommand{\Hom}{\mathsf{Hom}}
\newcommand{\eps}{\epsilon}

\title{Conditional upper bounds on the least character non-residue}
\author{Aritro Pathak}   
\date{}

\begin{document}
\maketitle

\begin{abstract}
   We extend known methods to establish upper bounds on the least character non-residues contingent on different zero-free regions within the critical strip, in particular on bounded rectangles within the critical strip along the $\sigma=1$ line at arbitrary heights. This relates to earlier conditional results on least character non-residues, and recent results of Granville and Soundararajan, on character sums.
\end{abstract}

\section{Introduction.}
\ \ \ \ \ For any non-principal character $\chi$ (mod $q$), consider  $n(\chi)$, the least character non-residue of $\chi$, which is the least positive integer $n$ such that $(n,q)=1$ and thus $\chi(n)\neq 0$, and further $\chi(n)\neq 1$.

For quadratic characters, Vinogradov's longstanding conjecture states that for any fixed $\eps>0$, and any prime modulus $p$, we must have $n(p)\ll p^{\eps}$. The corresponding conjecture on upper bounds on the character sums states that
\begin{align}\label{eq:eq1}
\sum\limits_{n\leq x} \chi(n)=o(x) \ \ \text{when} \ \ x\geq q^{\eps},
\end{align}
for some fixed $\eps$ and this implies Vinogradov's conjecture. Burgess's longstanding results\cite{Bur1,Bur2} give effective answers for $x\gg q^{1/4+\eps}$ when $q$ is cube-free and some weaker variants for general $q$(see, for example, equations 12.55, 12.56 and 12.57 in \cite{IK}).

Prior and recent works have focused on the connections between the zeros of $L$-functions and large character sums, or large character non-residues in different settings (see for example, \cite{Pier,Mei} and the references therein).

Assuming the Generalized Riemann Hypothesis(GRH), that all the zeros of $L(s,\chi)$ lie on the $\sigma=1/2$ line, one can immediately deduce a stronger upper bound of the form $n(\chi)\ll (\log q)^{2}(\log \log q)^{2}$ (see the discussion prior to Theorem 13.11 in \cite{Mont}). One can do even better and gain $n(\chi)\ll (\log q)^{2}$ under the same GRH hypothesis with the methods of \cite{Ankeny}(a simplification of this argument is presented in Theorem 13.11 of \cite{Mont}). Arguments in Rodosskii \cite{Rodo}, \cite{Mont2} Chapter 9, Theorem 1, show that weaker hypotheses on zero-free regions in the critical strip also give similarly strong upper bounds on the least character non-residue.
A similar result is shown in Theorem 10.6 of \cite{BD}.

As discussed in \cite{Sound}, observations of Heath-Brown on sequences of quadratic characters that violate the Burgess bound\cite{Bur0} in the form $q^{(1/4\sqrt{e}) +o(1)}$ would locate zeroes off the $\sigma=1/2$ line; further work of Banks and Makarov\cite{Banks} shows that if there is a sequence of quadratic characters for which the character sums are large in a smooth way, then one can locate zeros near $1$ of the corresponding $L$-functions.

In \cite{Sound}, it is shown that when one assumes that the primitive character sums modulo $q$ are large for $x\geq q^{\eps}$ in the sense that \cref{eq:eq1}fails, then a positive proportion of the zeros of the corresponding $L$-function lie off the $\sigma=1/2$ line and further the location of such zeros is approximately identified. There is also recent work by Tao \cite{Tao} relating the Vinogradov conjecture to the Elliot-Halberstam conjecture on the distribution of primes in arithmetic progressions.

We denote, as usual, the zeros $\rho$ of an $L$-function by $\beta+i\gamma$.  To motivate our first theorem, we recall a result of Backlund (Theorem 13.5 in \cite{Tit}), by which the Lindel\"{o}f hypothesis is equivalent to the statement that for every $\eps$, the number of zeros of the zeta function $\zeta(s)$ with real part at least $1/2 +\eps$ and imaginary part between $T$ and $T+1$ is at most $o(\log T)$. More precisely, we define,
\begin{align*}N(\eps,T):=|\{ \rho=\beta +i\gamma|\zeta(\rho)=0, \beta> \frac{1}{2} +\eps \}|.\end{align*}

Then the Lindel\"{o}f hypothesis is equivalent to the following:
\begin{theorem}
    For every $\eps>0$ we have, 
    \begin{align}
        N(\eps,T+1)-N(\eps, T)=o(\log T)
    \end{align}
\end{theorem}

A strong quantitative version of this hypothesis is that for some $0<\delta<1/2$, all the zeros of $\zeta(s)$ lie within the region $\mathcal{R}'=\mathcal{R}\cup \{\sigma+it| \delta\leq \sigma\leq 1-\delta, |t|<\theta \}$, where
\begin{equation}\label{eq;eq1}
\mathcal{R}=\{\sigma+it||\sigma -1/2|<\Phi(|t|), |t|\geq \theta\},
\end{equation}
 and $\Phi(t)$ is a function that tends to $0$ slowly as $t\to \infty$. 

Here in \cref{thm2}, for the case of the $L(s,\chi)$ functions, we have a weaker hypothesis than the above; with a quasi-Generalized Riemann Hypothesis condition. We have a zero-free region which is the infinite strip where $1-\delta< \sigma<1$ and its symmetric image about  the $\sigma=1/2$ line. 

We follow the idea of Ankeny's original proof (see for reference, Theorem 13.11 in \cite{Mont}), and it is a simpler argument than that used to prove the ensuing Theorem 2.
\begin{theorem}\label{thm2}
    Suppose that $\chi$ is a non principal character modulo $q$ and we have $L(s,\chi)\neq 0$ in $\mathcal{R}''(\delta)$, where $0 \leq \delta< \frac{1}{2}$
    \begin{align}
    \mathcal{R}'' (\delta):=\{\sigma+it: 1-\delta<\sigma <1 \}.
\end{align} Then $n(\chi)\ll (\log q)^{\frac{1}{\delta}}$.
\end{theorem}
 Using inverse Mellin transform techniques, we also prove the following stronger result, using a more restricted zero free region, generalizing the technique of an earlier result of Rodosskii\cite{Rodo}; see also \cite{Mont2} Chapter 9 Theorem 1 for an exposition.


\begin{theorem}\label{mainthm}

There exist absolute constants $C>0$ and $K_1>0$ such that the following holds.
Let $|t_0|>1$, let $\chi$ be a non-principal character $(\bmod\, q)$ with
$q \ge e^{|t_0|}+4$, and let
\[
\frac{1}{\log(|t_0|+4)} \le \delta \le \frac{1}{2}.
\]
If
\[
L(s,\chi)\neq 0 \quad \text{for } 1-\delta<\sigma<1, \qquad
|t| \le K_1 |t_0|^2 \log q,
\]
then
\[
n(\chi) < C\left( K_1 |t_0|^2 \frac{\log q}{\delta} \right)^{1/\delta}.
\]
\end{theorem}

\begin{remark}In \cref{mainthm} we will actually use the hypothesis that $L(s,\chi)\neq 0 \quad \text{for } 1-\delta<\sigma<1, 
|t-t_0| \le K_1 |t_0|^2 \log q$, which implies that the zero free region is of the form stated in \cref{mainthm} upto a change of constant.\end{remark}

This shows that choosing any $t_0\in \R, |t_0|> 1$, and sufficiently large $q$, one can pinpoint zeros of the $L(s,\chi)$ function in the rectangles of the form $\{\sigma +it| 1-\delta\leq  \sigma\leq 1, |t-t_0|\leq K_1|t_0|^{2}\log q \}$, if the least character non-residue for $\chi(\text{mod q})$ fails the bound established in this theorem. In the case of $t_0=0$, the result of Rodosskii establishes a  bound with a zero-free region of smaller height. Also, Theorem 1.3 of \cite{Sound} gives a bound on character sums $S(x,\chi):=\{\sum\limits_{n\leq x}\chi(n)\}$ in the case $x\approx \exp(\sqrt{\log q})$, when a rectangle in the critical strip of dimensions similar to the one in \cref{mainthm} does not contain too many zeroes.

Here, we also state the following two lemmas that would be used in several places in the proof of \cref{thm2,mainthm}, and which are stated as Theorem 10.17 and Theorem 11.5 of \cite{Mont}.

\begin{lemma}\label{lemma2}
Let $\chi$ be a character modulo $q$. The number of zeros 
$\rho = \beta + i\gamma$ of $L(s,\chi)$ in the rectangle
\[
0 \le \beta \le 1, \qquad T \le \gamma \le T+1
\]
is $\ll \log\bigl(q(|T|+2)\bigr)$.
\end{lemma}

\begin{lemma}\label{lemma}
    Let $n(r;t,\chi)$ denote the number of zeros $\rho$ of $L(s,\chi)$ in the disc 
$|\rho-(1+it)|\le r$. Then
\[
n(r;t,\chi)\ll r\log(q\tau)
\]
uniformly for
\[
\frac{1}{\log(q\tau)} \le r \le \frac{3}{4}.
\]
\end{lemma}

\section{Proofs.}

\begin{proof}[Proof of \cref{thm2}]

Fix an arbitrarily small $\delta>0$. We consider the expression, 
\begin{equation}\label{eq;eq2}
    \sum_{n\leq x} \chi(n)\Lambda(n)(x-n)=\frac{-1}{2\pi i}\int_{\sigma_0 -i\infty}^{\sigma_0 +i\infty} \frac{L'}{L}(s,\chi)\frac{x^{s+1}}{s(s+1)}ds, \ \ \sigma_{0}>1.
\end{equation}

\noindent We pull the contour to the line $\sigma=\delta/2$, say, and find a contribution from the zeros which is 
    \begin{align}\label{eq:estimate}
        -\sum\limits_{\rho}\frac{x^{\rho+1}}{\rho(\rho+1)},
    \end{align}
where, under the hypothesis, the sum is over all the non-trivial zeros of the $L$-function. Using a standard count for the number of zeros of the $L$ function as in \cref{lemma2}, this will be bounded in magnitude by 
\begin{align}\label{eqimp2}
\ll x^{2-\delta}\log q. 
\end{align}

Therefore we have, 

The line integral where $\sigma=\delta/2$, gives us,
\begin{equation}\label{eq:equnimp}
 -\frac{x^{1+ \delta/2}}{2\pi i}\int_{-\infty}^{\infty} \frac{L'}{L}\Big(\frac{\delta}{2} +it,\chi\Big) \frac{x^{it}}{(\delta/2 +it)(1+\delta/2 +it)}dt.
\end{equation}

\noindent This line integral will have a smaller contribution than that in \cref{eq:estimate}. This follows by considering Eq 12.5 in \cite{Mont} and noting that
\begin{align}
    \frac{L'}{L}\Big(\frac{\delta}{2} +it \Big)\ll \frac{1}{\delta}\log(q(|t|+4)),
\end{align}
and thus upon integrating we get from \cref{eq:equnimp}  a contribution that is bounded by,
\begin{align}\label{eqimp}
\ll \frac{1}{\delta}x^{1+\delta/2}\log q \end{align}
It will be enough later to require that 
\begin{align}\label{eq:eq12}
x^{2-\delta}\gg \frac{1}{\delta}x^{1+\delta/2}.
\end{align}
 \ \ \  \ \ \ \       On the other hand, by use of partial summation, standard upper bound estimates for $\omega(n)$ (see for example, Theorem 2.10 of \cite{Mont}) and elementary estimates for the von-Mangoldt function we have for $x\geq C(\log q)(\log \log q)$ with $C$ large enough. In this case, we get $\frac{1}{2}C^{2}(\log q)^2 (\log \log q)^2$ for the main term on the right of \cref{eq8} and a term $O(C (\log q)^2 (\log \log q)^2)$ for the error term which is negligible in comparison, when $C$ is large enough (See for example, equation 13.29 of \cite{Mont}).
\begin{equation}\label{eq8}
    \sum_{n\leq x} \chi_0 (n)\Lambda(n)(x-n)=\sum_{n\leq x} \Lambda(n)(x-n)+O(x(\log x)(\log q))\gg x^{2}.
\end{equation}

\noindent Thus, when we have $\text{max}(x^{2-\delta}, (1/\delta)x^{1+\delta/2})\log q=x^{2-\delta}\log q\ll x^{2}$, or in other words,
\begin{align}\label{eq9}
   x\gg(\log q)^{\frac{1}{\delta}},
\end{align}
and considering the contributions from \cref{eqimp,eqimp2}, $ \chi$ cannot equal $\chi_0$ till the value of $x$. Thus, we are forced to have $n(\chi)\ll (\log q)^{\frac{1}{\delta}}$, when we also satisfy \cref{eq:eq12}. For \cref{eq:eq12} to be satisfied, it is enough to require, using \cref{eq9}, that, 
\begin{align}
    (\log q)^{\frac{1}{\delta}(1 -\frac{3\delta}{2})}\gg \frac{1}{\delta}\Leftrightarrow \delta(\log q)^{\frac{1}{\delta}-\frac{3}{2}}\gg 1.
\end{align}

But this is true for all $q\geq 2$ and $\delta\leq \frac{1}{2}$.
\end{proof}

Next we prove \cref{mainthm}. 

\bigskip

\begin{proof}[Proof of \cref{mainthm}]
 Let $x(q,t_0)=\Big(\frac{K_2 t_{0}^{2}\log q}{\delta}\Big)^{\frac{1}{\delta}},y(q,t_0)=\sqrt{C}$  for constants $K_2,C$ , such that $2\leq y(q,t_0)\leq x(q,t_0)^{1/3}$. Henceforth, the dependence on $q, t_0$ of these $x,y$ parameters should be understood, and we don't write them as functions of $q,t_0$ explicitly.

 \bigskip
 
As will be seen in the course of the proof, the kernel function that works for our purpose, is
\begin{equation}
    K(s,t_0)=\int_{\frac{x}{y^{2}}}^{xy^{2}} \Big(2 \log y -\Big|\log \frac{x}{u}\Big|\Big)\Big(\frac{x}{u}\Big)^{1+it_0}u^{s-1} du,
\end{equation}
where the weight function is $w(u)=\big(2 \log y -|\log \frac{x}{u}|\big)\big(\frac{x}{u}\big)^{1+it_0}$


This weight function is  non-zero within the interval $(x/y^{2}, xy^{2})$ and take $w(u)=0$ elsewhere.

With a change of variable: $\frac{x}{u}=m$, this integral becomes 
\begin{multline}
    K(s,t_0)=x^{s}\int_{\frac{1}{y^{2}}}^{y^{2}} \Big( 2\log y -|\log m|  \Big) m^{it_0 -s} dm=2\log y\cdot x^{s} \Big( \frac{m^{it_0 -s+1}}{it_0 -s+1} \Big)|_{\frac{1}{y^{2}}}^{y^{2}}\\
    +x^{s}\int_{1/y^{2}}^{1} (\log m) m^{it_0 -s}dm -x^{s}\int_{1}^{y^{2}}(\log m) m^{it_0 -s} dm.
\end{multline}

\noindent This gives us:

\begin{multline}
    K(s,t_0)=\frac{2(\log y) x^{s} y^{2(it_0 -s +1)} }{it_0 -s+1}-\frac{2(\log y) x^{s} y^{-2(it_0 -s +1)} }{it_0 -s+1} -\frac{2x^{s}}{(it_0 -s+1)^{2}}\\ +\frac{2(\log y) x^{s} y^{-2(it_0 -s +1)} }{it_0 -s+1}-\frac{2(\log y) x^{s} y^{2(it_0 -s +1)} }{it_0 -s+1}+\frac{x^{s}y^{-2(it_0 -s+1)}}{(it_0-s+1)^{2}}+\frac{x^{s}y^{2(it_0 -s+1)}}{(it_0-s+1)^{2}}.
\end{multline}

Upon cancellations among four terms above, we get, 
\begin{align}
    K(s,t_0)=\frac{-2x^s +x^s y^{-2(it_0 -s+1)} +x^s y^{2(it_0 -s +1)} }{(it_0 -s+1)^2}.
\end{align}

\noindent Thus, we have for the kernel $K(s,t_0)$,

\begin{equation}
K(s,t_0)=x^{s}\Big(\frac{y^{s-it_0 -1}- y^{1+it_0 -s}}{s-it_0 -1}\Big)^{2}.
\end{equation}


\bigskip 

 Let $n(r;t,\chi)$ denote the number of zeros of $L(s,\chi)$ in the disk $D= \{x\in \mathbb{C}:|x-(1+it)|<r\}$. For any $t_0 \in \mathbb{R}$ with $\frac{1}{\log q(|t_0|+4)}\leq R \leq 1$, we have the estimate,
\begin{equation}\label{impeq}
    \sum\limits_{|\rho -1-it_0|>R} \frac{1}{|\rho- 1-it_0|^{2}} \ll \frac{\log q(|t_0|+4)}{R}.
\end{equation}

\noindent This follows from the fact that 

\begin{equation}
    \sum\limits_{R<|\rho -1-it_0|<2R} \frac{1}{|\rho -1-it_0|^{2}} \ll \frac{1}{R^{2}} n(2R;t_0,\chi) \ll \frac{\log q(|t_0|+4)}{R},
\end{equation}
where the first inequality follows because each individual summand on the left is bounded by $1/R^{2}$ and the number of terms is bounded by $n(2R;t_0,\chi)$. From the estimates on $n(R;t_0,\chi)$ from \cref{lemma}, 
we have the second inequality. Now replacing $R$ by $2^{k}R$, and summing over sufficiently many terms, we have a count of all the zeros within the annulus $R<|\rho -1-it_0|<1$, and we incorporate a constant pre-factor arising from the finite geometric progression, to get,
\begin{equation}\label{eq16}
    \sum\limits_{R<|\rho -1-it_0|<1} \frac{1}{
    |\rho -1-it_0|^{2}}\ll \frac{\log q(|t_0|+4)}{R}.
\end{equation}

For zeros further away from the point $(1+it_0)$, we use the estimate that the number of zeros of the function $L(s,\chi)$ within the range $0\leq \beta \leq 1$ and $T\leq \gamma \leq T+1$ is $\ll \log (q(|T|+4))$. Therefore putting $K=1$ in \cref{otherimp} below, it can be seen using \cref{eq16} along with an easy comparison that the bound in \cref{impeq} is attained. 

The zero-free region in the hypothesis of the theorem is symmetric with respect to the line $y=t_0$, and it is clear that for any positive integer $K$, we attain an upper bound of the form,

\begin{multline}\label{otherimp}
     \sum\limits_{|\rho-1-it_0|>K} \frac{1}{|\rho-1-it_0|^{2}} \ll \sum\limits_{n=K}^{\infty} \frac{\log q(n+|t_0|+4)}{n^{2}} \\ \ll \frac{\log q}{K}+\frac{\log(K+|t_0|+4)}{K}+ \Big(\frac{1}{|t_0|+4}\Big)\log\Big(1+ \frac{|t_0|+4}{K}\Big).
\end{multline}


The weight function being $w(u)= (2\log y -|\log(x/u)|)(x/u)^{1+it_0}$ in the interval $2\leq y \leq x^{1/3}$, we get for the weighted sum\footnote{See for example, the corresponding equation 13.31 in \cite{Mont}.},

\begin{equation}\label{eq:mainone}
\sum\limits_{n} w(n)\chi(n)\Lambda(n)= \frac{-1}{2\pi i} \int_{\sigma_0 -i\infty}^{\sigma_0 +i\infty} \frac{L'}{L}(s,\chi) \Bigg(\frac{y^{s-1-it_0}-y^{1+it_0-s}}{s-it_0-1} \Bigg)^{2} x^{s}ds.
\end{equation}

\noindent Here we have taken $\sigma_0 >1$ since $L(s,\chi)$ may have zeros in the critical strip and thus $\frac{L'(s,\chi)}{L(s,\chi)}$ would have simple poles in the critical strip. We move the contour to the left of the imaginary axis, to the line $\sigma_0 = -1+\theta$ for some small $\theta$ and get from \cref{eq:mainone},

\begin{multline}\label{eq:imp2}
    -\sum\limits_{\rho}\Big( \frac{y^{\rho -1-it_0}-y^{1+it_0-\rho}}{\rho -1-it_0}  \Big)^{2}x^{\rho} -\frac{(1-\kappa)}{(1+it_0)^{2}}(y^{1+it_0}- 1/y^{1+it_0})^{2} \\ -\frac{1}{2\pi i} \int_{-1+\theta -i\infty}^{-1+\theta +i\infty} \frac{L'}{L}(s,\chi) \Big( \frac{y^{s-1-it_0}-y^{1+it_0-s}}{s-it_0-1} \Big)^{2}x^{s}ds.
\end{multline}

\noindent The second term arises because $L(s,\chi)$ has a zero at $s=0$ when $\chi(-1)=1$ in which case $\kappa=0$ and otherwise $\kappa=1$.

If $\chi$ is induced by a primitive character $\chi^{*}$, then from the standard expression for $L(s,\chi)$ in terms of $L(s,\chi^{*})$, we take the logarithmic derivative and get that,

\begin{equation*}
    \frac{L'}{L}(s,\chi)=\frac{L'}{L}(s,\chi^{*})+\sum\limits_{p|q} \frac{\chi^{*}(p)\log p}{p^{s}-\chi^{*}(p)}.
\end{equation*}

In the sum above, we have $\sigma=-1+\theta$, and so it is easily seen that the summand is $\ll \log p$ and thus the sum of $\log p$ over all the primes $p|q$ is trivially bounded by $\log q$ itself. Further, since we have chosen the contour line to be where $\sigma_0=-1+\theta$, by another standard result (See for example Theorem 12.9 of \cite{Mont}), we can bound the above logarithmic derivative term for the primitive character by $\ll \log 2q(|t|+4)$, and then easily get a finite contribution from the contour integral and conclude that the total contribution to \cref{eq:imp2} from these terms is bounded by
\begin{equation*}
\ll x^{-1+\theta}y^{4-2\theta} \log q(|t_0|+4)\leq x^{-1+\theta}y^{4} \log q(|t_0|+4).
\end{equation*}

\noindent This will turn out to be smaller than the contribution from the first sum in \cref{eq:imp2}: recall that we started out with the bound $2\leq y \leq x^{1/3}$, and we will consider $y=O(1)$.  

Now first consider the contribution to the sum in \cref{eq:imp2} from the zeros $\rho=\beta +i\gamma$, with $\beta \leq 1-\delta$. In this case, considering $R=\delta$ in \cref{impeq}, we get the bound,

$$\ll \sum\limits_{|\rho-1-it_0|>\delta} \frac{x^{1-\delta}(3+y^{2})}{|\rho -1-it_0|^{2}} \ll x^{1-\delta}(3+y^{2}) \frac{(\log q(|t_0|+4))}{\delta}.$$

\noindent Further if $\rho$ is a zero in the critical strip with $\beta >1-\delta$, then by hypothesis, we have $|\gamma-t_0|\geq K(q,t_0)=K_{1}|t_0|^{2}\log q$. In order to use \cref{otherimp}, one can choose $K_1$ large enough so that $ K=\lfloor K(q,t_0) \rfloor$ is a positive integer. Each of these summands in \cref{eq:imp2} is bounded from above by $x(3+y^{2\delta})/(|\rho-1-it_0|^{2})$, when $1-\delta<\beta<1$.


We see that when $x$ is large enough, both of these bounds are greater than the bound of $x^{-1+\theta}y^{4}\log q(|t_0|+4)$ obtained from the contour integral in \cref{eq:imp2}. Combining with the estimate from \cref{otherimp}, we have that


\begin{multline}
     |\sum_{n}w(n)\chi(n)\Lambda(n)|\leq c\Bigg(\Big(x^{1-\delta}(3+y^{2}) \frac{(\log q(|t_0|+4))}{\delta}\Big)\\ + x(3+ y^{2\delta})\Bigg(\frac{\log q}{K}+ \Big(\frac{1}{|t_0|+4}\Big)\log\Big(1+ \frac{|t_0|+4}{K}\Big)+\frac{\log(K+|t_0|)}{K} \Bigg)\Bigg).
\end{multline}

\noindent If we now choose $x$ and  $K_1$ above so that $K=\delta x^{\delta}$ , then we will have $x(q,t_0)=\Big(\frac{K_2 t_{0}^{2}\log q}{\delta}\Big)^{\frac{1}{\delta}}$ for some constant $K_2$ as written in the beginning of the proof. We obtain,
\begin{multline}\label{eq:oneimp}
     |\sum_{n}w(n)\chi(n)\Lambda(n)|\leq c'\Big(x^{1-\delta}(3+y^{2}) \frac{(\log q(|t_0|+4))}{\delta}+\frac{1}{\delta}x^{1-\delta}(3+y^{2\delta})\big(\log (\delta x^{\delta}+|t_0|)+\log q \big)\\ +\Big(\frac{1}{|t_0|+4}\Big)\log\Big(1+ \frac{|t_0|+4}{\delta x^{\delta}}\Big) x(3+y^{2\delta})\Big).
\end{multline}

Next we estimate the sum  $\sum\limits_{n}w(n)\chi_0(n)\Lambda(n)$ where $\chi_0$ is the principal character modulo $q$. In that case, from \cref{eq:mainone}, putting in the principal character, we would then proceed as in the proof of the Prime Number Theorem, and then the leading order term would come from the residue of the function $-\frac{L'}{L}(s,\chi_0) \Big(\frac{y^{s-1-it_0}-y^{1+it_0-s}}{s-it_0-1} \Big)^{2}x^{s}$. We have a simple pole at $s=1$ and none at $s=1+it_0$. We have the residue at $s=1$ 

\begin{equation}\label{eq:last}
    -\frac{x}{t_{0}^{2}}\big(y^{-it_0} -y^{it_0}\big)^{2}=\frac{4x}{t_{0}^{2}}\sin^{2}(t_0 \log y).
\end{equation}
For the $y$ values we restrict to the sequence of values of $y_k$ so that $t_0 \log y_k =\frac{\pi}{2}(2k+1)$ (as well as $y_k\geq 2$) , i.e., we choose the values $y_k=e^{\frac{\pi}{2t_0}(2k+1)}$. Thus, $\sin^{2}(t_0 \log y_k)=1$ and then the leading order term in the sum $\sum\limits_{n}w(n)\chi_0(n)\Lambda(n)$ is $\frac{4x}{t_{0}^{2}}$.  Note that given any value of $t_0$ with $|t_0|>1$, we can choose $k$ appropriately so that $y(q,t_0)=O(1)$ and $y(q,t_0)\geq 2$.

We take $x$ large enough so that
\begin{equation}\label{eq24}
\delta x^{\delta}\gg (|t_0|+4),
\end{equation}
in which case the third term on the right of \cref{eq:oneimp} also has a factor of $x^{1-\delta}$, upon using the approximation of $\log\big(1+ \frac{|t_0|+4}{\delta x^{\delta}}\big)$by $ \frac{|t_0|+4}{\delta x^{\delta}} $. Given $|t_0|>1$, we consider $q\geq e^{|t_0|+4}$, in which case $\log q(|t_0|+4)$ is approximated by $ \log q$. In this case, noting that $y(q,t_0)=O(1)$, the dominant term on the right of \cref{eq:oneimp} is $c x^{1-\delta}\frac{\log q}{\delta}$.

If we require for $q\geq e^{|t_0|+4}$ , and with some large enough absolute constant $K_1$ that
\begin{align}
    \frac{4x}{t_{0}^{2}}\geq K_1 x^{1-\delta}\frac{\log q}{\delta}\Leftrightarrow 4x^{\delta}\delta \geq K_1 t_{0}^{2}\log q,
\end{align}
in which case \cref{eq24} is also satisfied, then the character $\chi(\text{mod} \ q)$ can't equal the principal character $\chi_0(\text{mod} \ q)$ till the value of $xy^{2}=Cx=C(\frac{K_2 t_{0}^{2}\log q}{\delta})^{\frac{1}{\delta}}$, due to the equality in Eq. 24 and the inequality in Eq. 23, and in which case we have $n(\chi)<C(\frac{K_2 t_{0}^{2}\log q}{\delta})^{\frac{1}{\delta}}$.
\end{proof}
\begin{remark}In this case we choose $\frac{1}{\log (|t_0|+4)}\leq \delta\leq \frac{1}{2}$, and thus for any integer $m$, where $2<m<\log (|t_0|+4)\leq \log \log q$, we have an upper bound on the least character non-residue in the form of an integral power of $\log q$, with a zero-free region of width $\frac{1}{m}$. 

Note that with a choice of $\frac{1}{|t_0|^{p}}<\delta<\frac{1}{2}$ , and thus with a possibly narrower zero-free region for any positive integer $p$, one can see that the upper bound on the least character non-residue actually approaches the classically known bounds, which is thus not optimal.
\end{remark}

\section{Further remarks, and future directions.} 
Using these techniques, one can in future work with the following zero-free regions:
\begin{enumerate}
\item \begin{itemize}
\item For $i\geq 1$, let $\delta_i\to 0$  be a decreasing  and possibly infinite sequence. Whenever $\beta>1-\delta_i$, we require $|t_0 -\gamma|>K_i$, where $K_i$ is an increasing sequence. 
\item Further, when $\frac{1}{2} \leq \beta\leq 1-\delta_1$ we impose no restriction on the $\gamma$ coordinate.
\end{itemize}
\bigskip
Recall that we use the notation $\beta +i\gamma$ for the zeros of an $L$-function. In this situation, when using the techniques of Theorem 2, the corresponding version of Eq. (21) will have several terms on the right with factors of the form $x^{1-\delta_i},\text{and} \ \frac{1}{K_i}$, and the argument will need to get modified accordingly. 

\item On the other hand, in the spirit of the result of Backlund mentioned in the Introduction, one can work with stronger zero-free regions where, 
\begin{itemize}
    \item $\delta_i\to 0$ is a decreasing sequence, $K_i$ is an increasing sequence, so that whenever $|t_0 -\gamma|>K_i$, we have $|\beta -\frac{1}{2}|<\delta_i$.
\end{itemize}
\end{enumerate}
\bigskip



\bigskip



\begin{thebibliography}{}
\footnotesize

\bibitem[Ank52]{Ankeny} Ankeny, Nesmith Cornett. "The least quadratic non residue." \textit{Annals of mathematics} (1952): 65-72.

\bibitem[Ban16]{Banks} Banks, William D., and Konstantin Makarov "Convolutions with probability distributions, zeros of L-functions, and the least quadratic nonresidue." \textit{Functiones et Approximatio Commentarii Mathematici} 55, no. 2 (2016): 243-280.

\bibitem[BD04]{BD} Bateman, Paul T., and Harold G. Diamond. \textit{Analytic number theory: an introductory course}. World Scientific, 2004.

\bibitem[Bur57]{Bur0} Burgess, David A. "The distribution of quadratic residues and non-residues." \textit{Mathematika} 4, no. 2 (1957): 106-112.

\bibitem[Bur62]{Bur1} Burgess, David A. "On character sums and L‐series." \textit{Proceedings of the London Mathematical Society} 12, no. 1 (1962): 193-206.

\bibitem[Bur63]{Bur2} Burgess, David A. "On character sums and L-series. II." \textit{Proceedings of the London Mathematical Society} 13, no. 1 (1963): 524-536.

\bibitem[Gra18]{Sound} Granville, Andrew, and Kannan Soundararajan. "Large character sums: Burgess's theorem and zeros of L-functions." \textit{Journal of the European Mathematical Society (EMS Publishing)} 20, no. 1 (2018):1-14.

\bibitem[Hea15]{Pier} Heath-Brown, David R., and Lillian B. Pierce.  ``Burgess bounds for short mixed character sums''. \textit{Journal of the London Mathematical Society} 91, no. 3 (2015): 693–708.

\bibitem[IK04]{IK} Iwaniec, Henryk, and Emmanuel Kowalski. \textit{Analytic Number Theory}. Vol. 53. Colloquium Publications, American Mathematical Soc., 2004.

\bibitem[M14]{Mei} Chang, Mei-Chu. "Short character sums for composite moduli." \textit{Journal d'Analyse Mathématique} 123, no. 1 (2014): 1-33.

\bibitem[Mon94]{Mont2} Montgomery, Hugh L. Ten lectures on the interface between analytic number theory and harmonic analysis. No. 84. CBMS Regional Conference Series in Mathematics. American Mathematical Society, 1994.

\bibitem[Mon07]{Mont} Montgomery, Hugh L., and Robert C. Vaughan. Multiplicative number theory I: Classical theory. No. 97. Cambridge Studies in Advanced Mathematics, Cambridge university press, 2007.


\bibitem[Rod56]{Rodo} Rodosskii, Kirill Andreevich. "On non-residues and zeros of L-functions." \textit{Izvestiya Rossiiskoi Akademii Nauk. Seriya Matematicheskaya} 20, no. 3 (1956): 303-306. 

\bibitem[Tao15]{Tao} Tao, Terence. "The Elliott–Halberstam conjecture implies the Vinogradov least quadratic nonresidue conjecture." \textit{Algebra and Number Theory} 9, no. 4 (2015): 1005-1034.


\bibitem[Tit86]{Tit}Titchmarsh, Edward Charles, and D. R. Heath-Brown. The theory of the Riemann zeta-function. Oxford university press, 1986.

\end{thebibliography}
\end{document}